\documentclass[reqno, 12pt]{amsart}
\usepackage{mathrsfs}
\usepackage{amssymb,latexsym}
\newtheorem{theorem}{Theorem}

\newtheorem{remark}[theorem]{Remark}
\newtheorem{proposition}[theorem]{Proposition}
\newtheorem{lemma}[theorem]{Lemma}
\newtheorem{definition}[theorem]{Definition}

\newtheorem{note}[theorem]{Note}

\newcommand{\ol}{\overline}

\newcommand{\abs}[1]{ |#1| }
\newcommand{\px}{ p(x) }
\newcommand{\norm}[1]{||#1||}

\begin{document}
\setlength{\baselineskip}{1.2\baselineskip}
\title  [High Energy Solutions to $p(x)-$Laplacian Equations of Schr\"{o}dinger type]
{High Energy Solutions to $p(x)-$Laplacian Equations of
Schr\"{o}dinger type}

\author{Duchao Liu$^*$}
\address{$^*$Department of Mathematics\\
        Lanzhou University, Lanzhou, 730000, China\\
        }
\email{liuduchao@yahoo.cn}

\author{Xiaoyan Wang$^\dag$}
\address{$^\divideontimes$Department of Mathematics\\
         Indiana University Bloomington\\
         IN, 47405, U.S.A}
\email{wang264@indiana.edu}
\thanks{$^{\dag, \bigstar}$Corresponding authors: X. Wang and J. Yao. }
\author{Jinghua Yao$^\bigstar$}
\address{$^\bigstar$Department of Mathematics\\
        Indiana University Bloomington\\
         IN, 47405, USA}
\email{yaoj@indiana.edu}

\maketitle

\begin{abstract}
In this paper, we investigate nonlinear Schr$\ddot{o}$dinger type
equations in $R^N$ under the framework of variable exponent spaces.
We propose new assumptions on the nonlinear term to yield bounded
Palais-Smale sequences and then prove the special sequences we find
converge to critical points respectively. The main arguments are
based on the geometry supplied by Fountain Theorem. Consequently, we
show that the equation has a sequence of solutions with high
energies.


\end{abstract}

\section{Introduction}

\setcounter{theorem}{0}
\setcounter{equation}{0}

Inspired by X.L. Fan \cite{12, 13} and L. Jeanjean \cite{24}, we
study the following nonlinear Schr$\ddot{o}$dinger type equation:
\begin{equation}\label{equ1}
-div( \abs{Du}^{p(x)-2}Du) + V(x)\abs{u}^{p(x)-2}u = f(x, u), x\in R^N,
\end{equation}
where $-div( \abs{Du}^{p(x)-2}Du)$ is called minus $p(x)$-Laplacian
and $V(x)$ satisfies the following condition $(V): V(x)\in C(R^N,
R), \inf_{x\in R^N}V(x)\geq V_0>0$; For every constant $M>0$, the
Lebesgue measure of the set $\{x\in R^N; V(x)\leq M\}$ is finite.
Here $V_0$ is a constant.

These equations involving the $p(x)$-Laplacian (also called
$p(x)$-Laplacian equations) arise in the modeling of
electrorheological fluids (see [2, 5, 31] and [27]) and image
restorations among other problems in physics and engineering. Lots
of classical equations, for example the classical fluid equations,
are also studied in this general framework (see the new monograph
[42] and the references therein). Different from the Laplacian
$\Delta:=\sum_j\partial^2_j$ (linear and homogeneous) and the
$p$-Laplacian $\Delta_p u(x):= div( \abs{Du}^{p-2}Du)$ (nonlinear
but homonegeous) where $0<p<\infty$ is a positive number, the
$p(x)$-Laplacian is nonlinear and nonhomogeneous. Besides the
applications we mentioned at the beginning of this paragraph, the
$p(x)$-Laplacian equations can be regarded as a nonlinear and
nonhomogeneous mathematical generalization of the stationary
Schr\"{o}dinger equation $\mathcal{H}u(x)=0$ where the Hamiltonian
is usually given by $\mathcal{H}:=-\frac{\hbar^2}{2m}\Delta+V(x)$.
For these connections and potential further generalizations, see
[45, 46, 47].

Next, we give the definitions of variable exponent spaces in order
to describe our problem precisely.

Let $\Omega$ be an open domain in $R^N$ and denote:
$$C_+(\Omega):= \{p(x)\in C(\Omega): 1<p^- :=\inf_{x\in \Omega}p(x) \leq p^+ := \sup_{x\in \Omega} p(x) < \infty \}.$$

For $p(x)\in C_+(\Omega)$, we consider the set:
$$ L^{p(x)}(\Omega) = \{ u: u \text{ is real-valued measurable function}, \int_\Omega \abs{u}^{p(x)}dx < \infty \}$$
We can introduce a norm on $L^{p(x)}(\Omega)$ by
$$\abs{u}_{p(x), \Omega} :=\inf\{k>0: \int_\Omega \abs{ \frac{u}{k} }^ {p(x)}dx\leq 1 \}$$
and $(L^{p(x)}(\Omega), \abs{\cdot}_{p(x), \Omega})$ is a Banach Space and we call it a variable exponent Lebesgue space.

Consequently, $W^{1, p(x)}(\Omega)$ is defined by
$$ W^{1, p(x)}(\Omega) = \{ u\in L^{p(x)}; \abs{Du} \in L^{p(x)}(\Omega)\}$$
with the following norm
$$\norm{u}_{p(x), \Omega} = \inf \{k>0; \int_\Omega \abs{\frac{Du}{k}}^{p(x)} + \abs{\frac{u}{k}}^{p(x)}dx\leq 1  \}.$$
Then $(W^{1, p(x)}{\Omega}, \norm{\cdot}_{p(x), \Omega})$ also
becomes a Banach space and we call it a variable exponent Sobolev
space.

Let
$$E:=\{u\in W^{1, p(x)}(R^N); \int_{R^N}V(x)\abs{u}^{p(x)}dx < \infty\}.$$
Then $E$ is a Banach space with the following norm
$$\norm{u}=\inf \{k>0; \int_{R^N} \abs{\frac{Du}{k}}^{\px} + V(x)\abs{\frac{u}{k}}^{\px} dx \leq 1\}.$$

Of course, our working space is $E$. Under reasonable and proper
assumptions, we shall show that $(\ref{equ1})$ has a sequence of
high energy solutions $\{u_n\}$ in $E$ in this paper (Theorem 2.2).

Since the last two decades, literatures and studies on variable
exponent spaces have sprung up like mushrooms (see [1, 2, 5, 7, 8,
9-20, 25, 31, 39-41]). These kinds of spaces are extensions of the
usual Lebesgue and Sobolev spaces
 $L^p(\Omega)$ and $W^{m, p}(\Omega)$ where $1\leq p< \infty$ is a constant. And they are special Orlicz spaces (see [26]). A lot of mathematical work has been done under the framework
 of the variable exponent spaces (see \cite{1, 4, 11, 27, 29, 35}). Meanwhile, a number of typical and interesting problems have
 come into light (see \cite{4, 6, 10, 15, 20, 22, 23, 28, 29, 33}). For example, local conditions on the exponent $p(x)$
 can assure the multiplicity of solutions to $p(x)-$Laplacian equations (See \cite{35}).

There is no doubt that there are mainly two characteristics when we
work with variable exponent spaces. For one thing, these spaces are
more complicated than the usual spaces (See \cite{3, 8, 17, 25}). As
a result, the related problems are more difficult to deal with. For
the other thing, we will obtain more general results if we work
under the framework of the variable exponent spaces because there
spaces are natural generalizations of the usual Sobolev and Lebesgue
spaces.

In X.L. Fan \cite{12}, the author considered a constrained
minimization problem involving $p(x)-$Laplacian in $R^N$. Under
periodic assumptions, the author could elaborately deal with this
unbounded problem by concentration-compactness principle of P. L.
Lions. In a following paper, X.L. Fan \cite{13}, the author
considered $p(x)-$Laplacian equations in $R^N$ with periodic data
and non-periodic perturbations. Under proper conditions, the author
was able to show the existence of solutions and gave a concise
description of the ground sate solutions. It is worth noting that
the periodicity assumptions are essential for the validity of
concentration-compactness principle under the framework of variable
exponent spaces (see the recent paper of Bonder and coworkers [43,
44] for the concentration-compactness theory in the variable
exponent space framework involving critical exponents). In our
paper, we also consider an unbounded problem. However, under
condition $(V)$, we could get some compact embedding theorems. In
fact, lots of other tricks can be used to recover some kinds of
compactness. For example, weight function method was used in
\cite{9}. In \cite{36}, we considered a combined effect of the
symmetry of the space and the coerciveness of potential $V(x)$.

We also want to mention the celebrated paper of L. Jeanjean
\cite{24}. In this paper, the author illustrated a completely new
idea to guarantee bounded $(PS)$ sequences for a given $C^1$
functional. Roughly speaking, we could consider a family of
functionals which contains the original one we are interested in.
When given additional structure assumptions, almost all the
functional in the family have bounded $(PS)$ sequences if the
functionals enjoy specific geometry properties! In fact, the
information of relevant functionals in the family can provide useful
information for the original functional. Under our conditions (see
Section 2), we could show the functional we consider satisfies the
fountain geometry. Then following L. Jeanjean's idea and Theorem 3.6
due to W. Zou \cite{38}, we could show equation ($\ref{equ1}$) has a
sequence of high energy solutions. We want to emphasize that our
condition $(C4)$ is somewhat mild and is first used in dealing with
$p(x)-$Laplacian equations. In addition, we do not need the
Ambrosetti-Rabinowits type condition here.

For the reader's convenience, we will recall some basic properties on the variable exponent spaces in the following part of this section.

\begin{proposition}(\cite{17, 18}) $L^{\px}(\Omega), W^{1,\px}(\Omega)$ are both separable, reflexive and uniformly convex Banach Spaces.
\end{proposition}

\begin{proposition}\label{Prop2}(\cite{17, 18}) Let $\rho(u)=\int_\Omega\abs{u(x)}^{\px} dx$ for $u\in L^{\px}(\Omega)$, then we have
    \begin{enumerate}
        \item $\abs{u}_{\px, \Omega}=1\Leftrightarrow \rho(u)=1;$
        \item   $\abs{u}_{\px, \Omega}\leq 1\Rightarrow \abs{u}^{p^+}_{\px, \Omega}\leq \rho(u)\leq \abs{u}^{p^-}_{\px, \Omega};$
        \item $\abs{u}_{\px, \Omega}\geq 1\Rightarrow \abs{u}^{p^-}_{\px, \Omega}\leq \rho(u)\leq \abs{u}^{p^+}_{\px, \Omega};$
        \item For $u_n\in L^{\px}(\Omega), \rho(u_n)\rightarrow 0 \Leftrightarrow \abs{u_n}_{\px, \Omega}\rightarrow 0$ as $n\rightarrow \infty;$
        \item For $u_n\in L^{\px}(\Omega), \rho(u_n)\rightarrow \infty \Leftrightarrow \abs{u_n}_{\px, \Omega}\rightarrow \infty$ as $n\rightarrow \infty.$
    \end{enumerate}
\end{proposition}

\begin{proposition}(\cite{17, 18, 26}) Let $\rho(u)=\int_\Omega \abs{Du(x)}^{\px} + \abs{u(x)}^{\px} dx$ for $u\in W^{1, \px}(\Omega)$, then we have
    \begin{enumerate}
        \item $\norm{u}_{\px, \Omega}=1\Leftrightarrow \rho(u)=1;$
        \item   $\norm{u}_{\px, \Omega}\leq 1\Rightarrow \norm{u}^{p^+}_{\px, \Omega}\leq \rho(u)\leq \norm{u}^{p^-}_{\px, \Omega};$
        \item $\norm{u}_{\px, \Omega}\geq 1\Rightarrow \norm{u}^{p^-}_{\px, \Omega}\leq \rho(u)\leq \norm{u}^{p^+}_{\px, \Omega};$
        \item For $u_n\in W^{1,\px}(\Omega), \rho(u_n)\rightarrow 0 \Leftrightarrow \norm{u_n}_{\px, \Omega}\rightarrow 0$ as $n\rightarrow \infty;$
        \item For $u_n\in W^{1,\px}(\Omega), \rho(u_n)\rightarrow \infty \Leftrightarrow \norm{u_n}_{\px, \Omega}\rightarrow \infty$ as $n\rightarrow \infty.$
    \end{enumerate}
\end{proposition}
The following property can be easily verified:

\begin{proposition} For $u\in E$, Let $\rho(u)=\int_{R^N}\abs{Du(x)}^{\px} + V(x)\abs{u(x)}^{\px} dx$. Then we have the following relations:
    \begin{enumerate}
        \item $\norm{u}=1 \Leftrightarrow \rho(u)=1;$
        \item $\norm{u}\leq 1 \Rightarrow \norm{u}^{p^+} \leq \rho(u) \leq \norm{u}^{p^-};$
        \item $\norm{u}\geq 1 \Rightarrow \norm{u}^{p^-} \leq \rho(u) \leq \norm{u}^{p^+}.$
    \end{enumerate}
\end{proposition}

From the above-mentioned properties, we can see that the norm and the integral (i.e. $\rho(u)$) don't enjoy the equality relation,
which is typical in variable exponent spaces and very different from the constant exponent case.

{\bf Notation.} For $\px\in C_+(\Omega), p^*(x)$ refers to the
critical exponent of $\px$ in the sense of Sobolev embedding,
that's, $p^*(x)=\frac{N\px}{N-\px}$ if $\px < N; p^*(x)=\infty$,
otherwise. For two continuous functions $a(x)$ and $b(x)$ in
$C(\Omega)$, $a(x)\ll b(x)$ means that $\inf_{x\in
\Omega}(b(x)-a(x)) > 0$. We will use the symbols
"$\rightharpoonup$", "$\rightarrow$" to represent weak convergence
and strong convergence in a Banach space respectively. While,
"$\hookrightarrow$", "$\hookrightarrow\hookrightarrow$" will be used
to denote continuous embedding and compact embedding between spaces
respectively. Technically, we use $C$ to denote a generic positive
constant.

\begin{proposition}(\cite{17, 18, 35})
    \begin{enumerate}
        \item Let $\Omega$ be a bounded domain in $R^N$. Assume that the boundary $\partial\Omega$ possesses cone property and $q(x)\in C(\ol{\Omega}, R)$ with $1\leq q(x) \ll p^*(x)$, then $W^{1, \px}(\Omega) \hookrightarrow \hookrightarrow L^{q(x)}(\Omega)$

        \item $W^{1, \px}(R^N) \hookrightarrow L^{q(x)}(R^N)$ if $p^+ < N$ and $q(x)\in C_+(R^N)$ satisfies $\px \leq q(x) \ll p^*(x)$.

    \end{enumerate}
\end{proposition}

Following the spirit of \cite{18}, we have the following proposition:

\begin{proposition}\label{Prop6} For $u\in E$, we define
$$ I(u)= \int_{R^N} \frac{1}{\px}( \abs{Du}^{\px} + V(x)\abs{u}^{\px})dx,$$
then $I\in C^1(E, R)$ and the derivative operator $L$ of $I$ is
$$\langle L(u), v\rangle=\int_{R^N}( \abs{Du}^{\px-2}Du\cdot Dv + V(x)\abs{u}^{\px-2}uv)dx, \forall u, v\in E$$
and we have:
\begin{enumerate}
    \item $L: E\rightarrow E^*$ (the dual space of $E$) is a continuous, bounded and strictly monotone operator;
    \item $L$ is a mapping of type ($S_+$), i.e. if $u_n \rightharpoonup u$ in $E$ and $\lim \sup_{n\rightarrow\infty}\langle L(u_n)-L(u), u_n-u\rangle\leq 0$, then $u_n\rightarrow u$ in $E$.
    \item $L: E\rightarrow E^*$ is a homeomorphism.
\end{enumerate}
\end{proposition}

\begin{proposition}\label{Prop7}(\cite{17, 18, 35}) Let $\Omega$ be a bounded domain in $R^N$. If $f(x, t)$ is a Caratheodory function and satisfies
$$\abs{f(x, t)} \leq a(x) + b\abs{t}^{ \frac{p_1(x)}{p_2(x)} }, \forall x\in \ol{\Omega}, t\in R$$
where $p_1(x), p_2(x)\in C_+(\Omega), b\geq 0$ is a constant, $0\leq a(x)\in L^{p_2(x)}(\Omega)$, then the superposition operator from $L^{p_1(x)}(\Omega)$ to $L^{p_2(x)}(\Omega)$ defined by $Su=f(x, u(x))$ is a continuous and bounded operator. Moreover, if $\Omega$ is unbounded (e.g. $\Omega=R^N$) and $a(x)\equiv 0$, the same conclusion is true.

\end{proposition}

In the variable Lebesgue space case, H$\ddot{o}$lder type inequality still holds.
\begin{proposition}\label{Prop8}(\cite{17}) Let $\Omega$ be a domain in $R^N$ (either bounded or unbounded) and $u\in L^{\px}(\Omega), v\in L^{p'(x)}(\Omega)$ where $p'(x) := \frac{\px}{\px - 1}$ is the conjugate exponent of $\px \in C_+(\Omega)$. Then the following H$\ddot{o}$lder type inequality holds
$$\int_\Omega\abs{uv}dx \leq (\frac{1}{p^-} + \frac{1}{p'^-})\abs{u}_{\px, \Omega} \abs{v}_{p'(x), \Omega}.$$
\end{proposition}

We still use this inequality in the following section (for example, in Lemma 2.7).

This paper is divided into 3 sections. For the readers' convenience,
we have recalled some basic properties of the variable exponent
spaces $W^{1, \px}(\Omega), L^{\px}(\Omega)$ in this section. In
section 2, we will state our assumptions on the nonlinear term and
our main result. Meanwhile, we prove some useful auxiliary results
in this section. In our opinion, these results themselves are
interesting and important when we study variable exponent problems.
In Sections 3, we are devoted to proving the main result.

\section{Main result}

\setcounter{theorem}{0}
\setcounter{equation}{0}

In this section, we will first specify our assumptions on the nonlinear term $f$ and give some comments on these assumptions. Then we state the main result.

We assume the following assumptions:

$(C1) f\in C(R^N \times R, R)$ satisfies
$$\abs { f(x, t) } \leq C( \abs{t}^{\px - 1} + \abs{t}^{q(x) - 1}), \forall t\in R, x\in R^N$$
$$f(x, t)t \geq 0, \text{for }    t\geq 0, x\in R^N$$
$$p(x) \leq q(x) \ll p^*(x),   \forall x\in R^N.$$

$(C2)$ There exists a constant $\mu > p^+$ such that
$$\liminf _{\abs{t} \to \infty} \frac{f(x, t)t}{\abs{t}^\mu} \geq C_0\text{ (a constant) }>0, \text{  uniformly for } x\in R^N.$$

$(C3)$
$$\limsup_{\abs{t} \to 0} \frac{f(x, t)t}{\abs{t}^{p^+}} = 0, \text{  uniformly for } x\in R^N.$$

$(C4)$ Let $F(x, t)=\int_0^t f(x ,s )ds$ and
$$G(x, t) := f(x, t)t - p^-F(x, t);$$
$$H(x, t) := f(x, t)t - p^+F(x,t).$$

We assume $G$ and $H$ satisfy the monotonicity condition: there
exist two positive constants $D_1$ and $D_2$ such that
$$G(x ,t) \leq D_1G(x, s) \leq D_2H(x, s), \text{for } 0\leq t \leq s.$$

$(C5)$
$$f(x, -t) = -f(x, t), \forall t\in R, x\in R^N.$$

\begin{definition} We say $u\in E$ is a solution to the equation ($\ref{equ1}$) if for any $v\in E$, the following equality holds
$$\int_{R^N}\abs{Du}^{\px - 2 }DuDv + V(x)\abs{u}^{\px-2}uv dx = \int_{R^N} f(x, u)v dx.$$
\end{definition}

Define a functional $\Phi$ from $E$ to $R$:
$$\Phi(u) = \int_{R^N} \frac{1}{\px} ( \abs{Du}^{\px} + V(x)\abs{u}^{\px})dx - \int_{R^N}F(x, u)dx.$$
Under our assumptions, we know the functional is $C^1$ (Proposition 1.6, Lemma 2.7 in the following) and for $v\in E$,
$$\Phi'(u)v=\int_{R^N} \abs{Du}^{\px-2}DuDv + V(x)\abs{u}^{\px-2}uvdx - \int_{R^N}f(x, u)vdx.$$
So the critical points of the functional $\Phi$ are corresponding to the solutions of the equation ($\ref{equ1}$).

Next are some comments and analysis on the assumptions we give.
\vspace{10pt}

1. $(C1)-(C4)$ are compatible. We give two examples. Let $f(x, t) = \abs{t}^{q(x) - 2}t$ with $q(x)\in C_+(R^N)$
satisfying $q(x)\ll p^*(x), q_->p^+$. Obviously, $(C1), (C2), (C3), (C5)$ hold. In order to verify $(C4)$,
we know that $F(x, t)=\frac{ \abs{t}^{q(x)}}{q(x)}, f(x, t)t=\abs{t}^{q(x)}$. Consequently,
$G(x, t)=(1- \frac{ p^-}{q(x)}) \abs{t}^{q(x)}, H(x, t)= (1-\frac{p^+}{q(x)})\abs{t}^{q(x)}$.
It's easy to verify that $G(x, t)$ is nondecreasing in $t\geq 0$. Therefore, $G(x, t)\leq G(x, s)$ if $0\leq t\leq s$. In view of $G, H\geq 0$, we know that
$$\frac{G(x, s)}{H(x, s)} = \frac{ q(x)-p^-}{ q(x)- p^+} \leq \frac{ q^+ - p^-}{ q^- - p^+ }.$$
Choosing $D_2=  \frac{ q^+ - p^-}{ q^- - p^+ }$, we obtain $G(x,
s)\leq D_2H(x, s)$ when $s\geq 0$. On the whole, $(C4)$ holds.

Next, we will illustrate another example. Let $f(x, t) =
\abs{t}^{q(x)-2}t \ln^a(\abs{t} + 1)$ where $q(x)$ satisfies
$q(x)\ll p^{*}(x)$, $q^{-}>p^{+}$ and $\epsilon>a>0$ is a real
number. In view of the following relations:
$$\lim_{\abs{t} \to \infty} \frac{ \ln^a(\abs{t}+1)}{ \abs{t}^\epsilon} = 0 \ \ \    \forall a\geq 0, \epsilon>0;$$
$$\lim_{\abs{t} \to 0}\frac{ \ln^a(\abs{t}+1)}{ \abs{t}^\epsilon} = \infty  \ \ \   \forall a\geq 0, \epsilon>0.$$
we can verify $(C4)$ similarly. Evidently, $(C1), (C2), (C3), (C5)$ hold.

From the two examples, we know lots of functions satisfy our assumptions. As a result, our main result is somewhat general.
\vspace{10pt}

2. Condition $(C1)$ means that $f(x, t)$ is subcritical in the variable sense. Different from things in constant case (i.e. $p^+=p^-$), we need $q(x)\ll p^*(x)$.
\vspace{10pt}

3. Condition $(C4)$ is crucial for our proof. It's because of this
condition that we could obtain bounded Palais-Smale sequence
(bounded $(PS)$ sequences for short). We give this condition to $f$
other than the famous Ambrosetti-Rabinowitz type condition. However,
we could still get bounded $(PS)$ sequences via an indirect method.
Lots of authors have tried to weaken the Ambrosetti-Rabinowits type
condition and they can only get weak type $(PS)$ sequences (usual
the Cerami Condition). It's known that $(C5)$ is much weaker than
the Ambrosetti-Rabinowitz type condition in the constant exponent
case ($p^+=p^-$) (see \cite{21}). \vspace{10pt}

4. Condition $(C5)$ assures the functional $\Phi$ we defined before is an even functional. So the condition is necessary for us to take advantage of the fountain geometry.

In this paper, we assume condition $(V)$ always holds and $p^+ < N$.
Hence, we know $E\hookrightarrow W^{1, \px}(R^N).$ Consequently,
$E\hookrightarrow L^{\px} (R^N), E\hookrightarrow L^{q(x)}(R^N)$ if
$q(x)\in C_+(R^N)$ satisfies $\px \leq q(x) \ll p^*(x)$.
\vspace{10pt}

Now we can state our main result clearly.
\begin{theorem}\label{theo1} Under condition $(V)$ and $(C1)-(C5)$, the equation ($\ref{equ1}$) has a
sequence of solutions $\{u_n\}$. Moreover, the solutions we get have
high energies, i.e. $\Phi(u_n)\rightarrow \infty$ as
$n\rightarrow\infty$.
\end{theorem}

In order to make the exposition more concise , we will give some
auxiliary results, some of which are very useful themselves.
\begin{lemma}\label{lemma3} Let $\Omega$ be a nonempty domain in $R^N$ which can be bounded or unbounded. We also allow $\Omega=R^N$. Then $L^{\px}(\Omega) \cap L^{q(x)}(\Omega) \subset L^{a(x)}(\Omega)$ if $p(x), q(x), a(x)\in C_+(\Omega)$ and $p(x)\leq a(x) \leq q(x)$. Moreover, if $p(x) \ll a(x) \ll q(x)$, the following interpolation inequality holds for $u\in L^{\px}(\Omega) \cap L^{q(x)}(\Omega):$
\begin{equation}
\int_\Omega \abs{u}^{a(x)} dx \leq 2|{\abs{u}^{a_1(x)}}|_{m(x),
\Omega} |{\abs{u}^{a_2(x)}}|_{m'(x), \Omega},
\end{equation}
where
$$a_1(x)=\frac{ p(x)(q(x)-a(x))}{q(x)-p(x)}, a_2(x)=\frac{q(x)(a(x)-p(x))}{q(x)-p(x)};$$
$$m(x)=\frac{q(x)-p(x)}{q(x)-a(x)},  m'(x)=\frac{q(x)-p(x)}{a(x)-p(x)}.$$
\end{lemma}

{\bf Sketch of Proof.}  For $L^{\px}(\Omega) \cap L^{q(x)}(\Omega)$,
we have
$$\int_\Omega \abs{u}^{\px} dx < \infty,  \int_\Omega \abs{u}^{q(x)} dx < \infty.$$
Obviously, $\abs{u(x)}^{a(x)} \leq \abs{u(x)}^{\px} + \abs{u(x)}^{q(x)}$ for $x\in\Omega$. Hence, $\int_\Omega \abs{u}^{a(x)} \leq \int_\Omega \abs{u}^{\px} dx + \int_\Omega \abs{u}^{q(x)} dx  < \infty$, which means $u\in L^{a(x)}(\Omega)$. For the interpolation inequality, the readers can see \cite{20}. \qed

\begin{lemma}\label{lemma4}
Under the condition $(V)$, $E\hookrightarrow\hookrightarrow L^{\px}(R^N).$
\end{lemma}
{\bf Proof}. We have known that $E\hookrightarrow L^{\px}(R^N)$. Next, we assume $u_n\rightharpoonup 0$ in $E$. We need to show $u_n \rightarrow 0$
in $L^{\px}(R^N)$ to complete the proof. by Proposition 1.2, it suffices to verify $\int_{R^N} \abs{u_n}^{\px} dx \rightarrow 0$ as $n\rightarrow \infty$.
For any given $R>0$, we write
$$I(n):= \int_{R^N} \abs{u_n}^{\px} dx = \int_{B(0, R)}\abs{u_n}^{\px} dx + \int_{R^{N}\backslash B(0, R)}\abs{u_n}^{\px} dx := I_1(n) + I_2(n).$$
Since $E\hookrightarrow W^{1, \px}(R^N)$ and $W^{1, p(x)}(B(0,
R))\hookrightarrow\hookrightarrow L^{\px}(B(0, R))$, we have
$I_1(n)\to 0$ as $n\to \infty$.

For any constant $M>0$, Let $A=\{ x\in R^N\backslash B(0, R);
V(x)>M\}$ and $B=\{ x\in R^N \backslash B(0, R); V(x)\leq M\}$. Then
$\int_A \abs{u_n}^{\px} dx \leq \int_A
\frac{V(x)}{M}\abs{u_n}^{\px}dx\leq
\frac{1}{M}\int_{R^N}V(x)\abs{u_n}^{\px}dx\leq \frac{C}{M}$; Since
for the constant $M>0, mes\{x\in R^N; V(x)\leq M\}$ is finite, we
can choose $R>0$ large enough such that $mes\{x\in R^N \backslash
B(0, R); V(x)\leq M\} \to 0$. Consequently, $\int_B \abs{u_n}^{\px}
\to 0$.

Now Let $M\to \infty$ and $R\to\infty$, we have $I(n)\to 0$ as $n\to
\infty$. \qed

\begin{lemma}\label{lemma5}
Under the condition $(V)$, $E\hookrightarrow\hookrightarrow L^{a(x)}(R^N)$ if $a(x)\in C_+(R^N)$ and $p(x) \leq a(x) \ll p^*(x)$.
\end{lemma}
{\bf Proof.} Let $u_n\rightharpoonup 0$ in $E$. We need to show $u_n\rightarrow 0$ in $L^{a(x)}(R^N)$ to finish the proof.

First, we assume that $p(x)\ll a(x)\ll p^*(x)$. We can choose $q(x)\in C_+(R^N)$ such that $a(x)\ll q(x)\ll p^*(x)$. It's obvious that $E\hookrightarrow L^{q(x)}(R^N)$. In view of $p(x)\ll a(x) \ll q(x)$, we use Lemma $\ref{lemma3}$ with $\Omega=R^N$ and obtain
\begin{equation}\label{equ2}
\int_\Omega \abs{u_n}^{a(x)} dx \leq 2|{\abs{u_n}^{a_1(x)}}|_{m(x),
\Omega} |{\abs{u_n}^{a_2(x)}}|_{m'(x), \Omega},
\end{equation}
where the symbols are the same as those of Lemma $\ref{lemma3}$.

Let $\lambda_n := ||{u_n}|^{a_1(x)}|_{m(x), \Omega}$ and
$\mu_n:=||{u_n}|^{a_2(x)}|_{m'(x), \Omega}$. By Proposition
$\ref{Prop2}$, we have
$$ \int_{R^N} \abs{ \frac{\abs{u_n}^{a_1(x)}}{\lambda_n}}^{m(x)}dx = \int_{R^N} \frac{\abs{u_n}^{\px}}{\lambda_n^{m(x)}} dx = 1;$$
$$ \int_{R^N} \abs{ \frac{\abs{u_n}^{a_2(x)}}{\mu_n}}^{m'(x)}dx = \int_{R^N} \frac{\abs{u_n}^{q(x)}}{\mu_n^{m'(x)}} dx = 1.$$
From the two equalities and Lemma $\ref{lemma4}$, we know
$$\min\{ \lambda_n^{m^+}, \lambda_n^{m^-}\} \leq \int_{R^N}\abs{u_n}^{\px}dx \to 0,$$
$$\min\{ \mu_n^{m'^+}, \mu_n^{m'^-}\} \leq \int_{R^N}\abs{u_n}^{q(x)}dx \leq C.$$
Anyway, we have $\lambda_n\to 0$ as $n\to \infty$ and $0\leq \mu_n
\leq C$. So $(2.2)$ yields that $\int_{R^N} \abs{u_n}^{a(x)}dx \to
0$ as $n\to\infty$.

Next, we assume $p(x)\leq a(x)\ll p^*(x)$. We can choose $q(x)\in C_+(R^N)$ such that $a(x)\ll q(x)\ll p^*(x)$. By the arguments above, we have
$$\int_{R^N}\abs{u_n}^{q(x)}dx \to 0.$$
By Lemma $\ref{lemma3}$, Lemma $\ref{lemma4}$, we have
$$\int_{R^N}\abs{u_n}^{a(x)}dx \leq \int_{R^N}\abs{u_n}^{p(x)}dx + \int_{R^N}\abs{u_n}^{q(x)}dx \to 0.$$\qed

The following lemma can be considered as an extension of the result in M. Willem [34, Appendix A].
\begin{lemma}\label{lemma6}
Assume $1\leq p_1(x), p_2(x), q_1(x), q_2(x) \in C(\Omega)$. Let $f(x, t)$ be a Carath$\acute{e}$odory function on $\Omega\times R$ and satisfy
$$\abs{ f(x, t) } \leq a\abs{t}^{ \frac{p_1(x)}{q_1(x)}} + b\abs{t}^{ \frac{p_2(x)}{q_2(x)} },  (x, t)\in \Omega\times R,$$
where $a,b>0$ and $\Omega$ is either bounded or unbounded. Define a Carath$\acute{e}$odory operator by
$$Bu :=f(x, u(x)),  u\in \mathscr{H}:=L^{p_1(x)}(\Omega)\cap L^{p_2(x)}(\Omega)$$
Define the space $\mathscr{E} :=L^{q_1(x)}(\Omega) +
L^{q_2(x)}(\Omega)$ with a norm
$$\norm{u}_\mathscr{E} =\inf\{ |{v}|_{q_1(x), \Omega} + |{w}|_{q_2(x), \Omega} : u=v+w, v\in L^{q_1(x)}(\Omega), w\in L^{q_2(x)}(\Omega)\}.$$
If $\frac{ p_1(x)}{q_1(x)} \leq \frac{p_2(x)}{q_2(x)}$ for
$x\in\Omega$, then $B=B_1 + B_2$, where $B_i$ is a bounded and
continuous mapping from $L^{p_i(x)}(\Omega)$ to $L^{q_i(x)}(\Omega),
i=1,2$. In particular, $B$ is a bounded continuous mapping from
$\mathcal{H}$ to $\mathscr{E}$.
\end{lemma}
{\bf Proof.} Let $\psi: R\rightarrow [0, 1]$ be a smooth function
such that $\psi(t)=1$ for $t\in (-1, 1); \psi(t)=0$ for $t \notin
(-2, 2)$. Let
$$g(x, t)=\psi(t)f(x, t),  h(x, t)=(1-\psi(t))f(x, t).$$
Because $\frac{ p_1(x)}{q_1(x)} \leq \frac{ p_2(x)}{q_2(x)}$ for
$x\in\Omega$, there are two constants $d>0, m>0$ such that
$$\abs{g(x, t)} \leq d\abs{t}^{ \frac{ p_1(x)}{q_1(x)} }, \abs{h(x, t)} \leq m\abs{t}^{ \frac{ p_2(x)}{q_2(x)}}.$$
Define
$$B_1u=g(x, u), u\in L^{p_1(x)}(\Omega);  B_2u=h(x, u), u\in L^{p_2(x)}(\Omega).$$
Then by Proposition $\ref{Prop7}$, $B_i$ is a bounded and continuous
mapping from $L^{p_i(x)}(\Omega)$ to $L^{q_i(x)}(\Omega), i=1,2$.
It's readily to see that $B:=B_1+B_2$ is a bounded continuous
mapping from $\mathcal {H}$ to $\mathscr{E}$. \qed

From Lemma $\ref{lemma4}$, Lemma $\ref{lemma5}$, we know that the
condition $(V)$ plays an important role. In enables $E$ to be
compactly embedded into $L^{\px}(R^N)$ type spaces. Using Lemma
$\ref{lemma5}$ and Lemma $\ref{lemma6}$, we can prove the following

Lemma $\ref{lemma7}$.
\begin{lemma}\label{lemma7}
Under assumption $(V)$ and $(C1)$, the functional $J(u) = \int_{R^N}
F(x, u)\,dx$ on $E$ is a $C^1$ functional. Moreover, $J'$ is
compact.
\end{lemma}
{\bf Proof.} The verification that $J$ is a $C^1$ functional is
routine and we omit it here. We only show that $J'$ is compact.
Because $E\hookrightarrow\hookrightarrow L^{\px}(R^N)$ (Lemma
$\ref{lemma4}$) and $E\hookrightarrow\hookrightarrow L^{q(x)}(R^N)$
(Lemma $\ref{lemma5}$), any bounded sequence $\{u_k\}$ in $E$ has a
renamed subsequence denoted by $\{u_k\}$ which converges to $u_0$ in
$L^{\px}(R^N)$ and $L^{q(x)}(R^N)$. Using Lemma $\ref{lemma6}$ with
$p_1(x)=\px, q_1(x)=\frac{\px}{\px -1}, p_2(x)=q(x),
q_2(x)=\frac{q(x)}{q(x)-1}$ and $\Omega=R^N$, we have
$J'(u)v=\int_{R^N}(B_1u + B_2u)vdx$ for $v\in E$. Hence,
$B_1(u_k)\to B_1(u_0)$ in $L^{q_1(x)}(\Omega)$ and $B_2(u_k)\to
B_2(u_0)$ in $L^{q_2(x)}(\Omega)$. Then H$\ddot{o}$lder type
inequality (Propsition $\ref{Prop8}$) and Sobolev embedding (Lemma
$\ref{lemma5}$) assure $J'(u_k)\to J'(u_0)$ in $E^*$, i.e. $J'$ is
compact. This proves the Lemma. \qed \vspace{12pt}

For conveience, we give the definition of $(PS)_c$ sequence for $c\in R$.

\begin{definition}
Let $\Pi$ be a $C^1$ functional defined on a real Banach space $X$. Any sequence $\{u_n\}$ satisfying $\Pi(u_n)\to c$ and $\Pi'(u_n)\to 0$ is called a $(PS)_c$ sequence. In addition, we call $c$ here a prospective critical level of $\Pi$.
\end{definition}

\begin{remark}\label{remark1}
(See also \cite{14}) Under the assumption of Theorem $\ref{theo1}$,
we have the following comments. $\Phi(u)=I(u) + J(u)$ and $\Phi'(u)
= I'(u) + J'(u)$ for $u\in E$. Since $I'$ is of type $(S_+)$
(Proposition $\ref{Prop6}$) and $J'$ is a compact (Lemma
$\ref{lemma7}$), we can easily derive that $\Phi'$ is of type
$(S_+)$. It's well-known that any bounded $(PS)_c$ sequence of a
functional whose Fr$\acute{e}$chet derivative is of type $(S_+)$ in
a reflexive Banach space has a convergent subsequence and so does
$\Phi$ here.
\end{remark}

\section{Proof of Theorem $\ref{theo1}$}
\setcounter{theorem}{0}
\setcounter{equation}{0}

We will first state the Fountain Theorem before our proof.

Let $X$ be a Banach space with the norm $\norm{\cdot}$ and let
$\{X_j\}$ be a sequence of subspaces of $X$ with $\dim{X_j} <
\infty$ for each $j\in \mathbb{N}$. Further, $X=\ol{
\bigoplus_{j=1}^\infty X_j}, W_k := \bigoplus_{j=1}^k X_j, Z_k
:=\ol{ \bigoplus_{j=k}^\infty X_j}$. Moreover, for $k\in \mathbb{N}$
and $\rho_k > r_k>0$, we denote:
$$B_k=\{u\in W_k: \norm{u}\leq \rho_k\};  S_k=\{u\in Z_k: \norm{u}=r_k\};$$
$$ c_{k}:=\inf_{\gamma\in \Gamma_k}\max_{u\in
B_k}\Phi(\gamma(u)),\,\mbox{where}$$
$$\Gamma_k :=\{ \gamma \in C(B_k, X): \gamma \text{ is odd and } \gamma|_{\partial B_k} = id\}.$$

\begin{theorem}\label{theo2}
(\cite{34}, Fountain Theorem, Bartsch, 1992) Under the aforementioned assumptions, let $\Phi\in C^1(X, R)$ be an even functional. If for $k>0$ large enough, there exists $\rho_k > r_k > 0$ such that
\begin{equation}
(A)\ a_k :=\max\{ \Phi(u) : u\in W_k, \norm{u} = \rho_k\} \leq 0,
\end{equation}
\begin{equation}
(B)\ b_k :=\inf\{ \Phi(u): u\in Z_k, \norm{u} = r_k\} \to \infty \text{ as }k\to \infty.
\end{equation}
then $\Phi$ has a $(PS)_{c_k}$ sequence for each prospective critical value $c_k$ and $c_k\to \infty$ as $k\to\infty$.
\end{theorem}

\begin{definition}
Let $X$ be a Banach space, $\Phi\in C^1(X, R)$ and $c\in R$. The function $\Phi$ satisfies the $(PS)_c$ condition if any sequence $\{u_k\}\subset X$ such that
    \begin{equation}
        \Phi(u_n)\to c, \Phi'(u_n)\to 0
    \end{equation}
has a convergent subsequence.
\end{definition}

\begin{remark}
In fact, if the following condition $(C)$ holds
    \begin{equation}
     (C)\ \Phi \text{ satisfies the }(PS)_c \text{ condition for every } c>0,
    \end{equation}
the sequence $\{c_k\}$ in Theorem $\ref{theo2}$ is a sequence of
unbounded critical values of $\Phi$. However, the condition $(C)$
isn't necessary to guarantee $c_k$ is a critical level. We just need
$(PS)_{c_k}$ condition.
\end{remark}

In order to use the decomposition technique, we need a theorem on the structure of a reflexive and separable Banach space.

\begin{lemma}\label{lemma34}
(See [37, Section 17]) Let $X$ be a reflexive and separable Banach space, then there are $\{e_n\}_{n=1}^\infty \subset X$ and $\{f_n\}_{n=1}^\infty \subset X^*$ such that:
$$f_n(e_m) = \delta_{n,m} =
\begin{cases}
1, \text{    if } n=m \\
0, \text{    if } n\neq m
\end{cases}
$$

$X= \ol{ span } \{e_n : n=1, 2, \cdots , \}, X^*=\ol{span}^{W^*} \{ f_n : n=1, 2, \cdots, \}.$
\end{lemma}

For $k=1, 2, \cdots$, and $X=E$, we will choose:
$$X_j =\text{span} \{e_j\}, W_k=\oplus_{j=1}^k X_j,  Z_k=\ol { \oplus_{j=k}^\infty
X_j}.$$\\

In the following, we identify the Banach space $E$ and the functional $\Phi$ as those we consider.
 Next, we will prove the main result step by step. First, we give a useful lemma. For simplicity,
 we write $|{u}|_{\px, R^N}$ as $|{u}|_{\px}$ when $\Omega=R^N$ for $\px \in C_+(R^N)$.

\begin{lemma}\label{lemma35}
Let $q(x)\in C_+(R^N)$ with $\px \leq q(x) \ll p^*(x)$ and denote
    \begin{equation}
    \alpha_k=\sup\{ \abs{u}_{q(x)}: \norm{u} = 1, u\in Z_k\},
    \end{equation}
then $\alpha_k\to 0$ as $k\to \infty$.
\end{lemma}
{\bf Proof.} Obviously, $\alpha_k$ is decreasing as $k\to\infty$.
Noting that $\alpha_k\geq 0$, we may assume that $\alpha_k\to \alpha
\geq 0$. For every $k>0$, there exists $u_k\in Z_k$ such that
$\norm{u_k}=1$ and $\abs{u_k}_{q(x)}>\frac{\alpha_k}{2}$. By
definition of $Z_k$, $u_k\rightharpoonup 0$ in $E$. Then Lemma
$\ref{lemma5}$ implies that $u_k\to 0$ in $L^{q(x)}(R^N)$. Thus we
have proved that $\alpha=0$. \qed\\

Using lemma $\ref{lemma35}$, we can prove the following Lemma $\ref{lemma36}$:
\begin{lemma}\label{lemma36}
Under the assumptions of Theorem $\ref{theo2}$, the geometry conditions of the Fountain Theorem hold, i.e. $(A)$ and $(B)$ hold.
\end{lemma}
{\bf Proof.} By $(C2)$ and $(C3)$, for any $\epsilon > 0$, there exists a $C(\epsilon) > 0$ such that
$$f(x, u)u\geq C(\epsilon)\abs{u}^\mu - \epsilon\abs{u}^{p^+}.$$
In view of $(C5)$, we have a constant, still denoted by $C(\epsilon)$, such that
$$F(x, u) \geq C(\epsilon)\abs{u}^\mu - \epsilon\abs{u}^{p^+}.$$
When $\norm{u} > 1$, we have
\begin{equation}
\begin{array}{rl}
\Phi(u) &=\int_{R^N} \frac{1}{\px} (\abs{Du}^{\px} + V(x)\abs{u}^{\px}) dx - \int_{R^N}F(x, u)dx
    \\&\leq \frac{1}{p^-}\norm{u}^{p^+} - C(\epsilon)\int_{R^N}\abs{u}^{\mu}dx + \epsilon\int_{R^N} \abs{u}^{p^+}dx.
\end{array}
\end{equation}
Let $u\in W_k$, since $\dim(W_k) < \infty$. all norms on $W_k$ are
equivalent. Hence $\Phi(u) \leq C\norm{u}^{p^+} - C\norm{u}^\mu$.
Because $\mu > p^+$, we can choose $\rho_k > 0$ large enough such
that $\Phi(u)\leq 0$ when $\norm{u}=\rho_k$. We have shown $(A)$
holds.

To verify $(B)$, we can still let $\norm{u} > 1$ without loss of generality. By $(C1)$ and $(C3)$, for any $\epsilon > 0$, there exists a $C=C(\epsilon) > 0$ such that
$$\abs{F(x, u)} \leq \epsilon\abs{u}^{p^+} + C\abs{u}^{q(x)},$$
So
\begin{equation}
\begin{array}{rl}
\Phi(u) &= \int_{R^N} \frac{1}{p(x)}( \abs{Du}^{\px} + V(x)\abs{u}^{\px}dx ) - \int_{R^N} F(x, u)dx
    \\& \geq \frac{1}{p^+} \norm{u}^{p^-} - \epsilon \abs{u}_{p^+}^{p^+} - C\max\{ \abs{u}_{q(x)}^{q^-}, \abs{u}_{q(x)}^{q^+}\}.
 \end{array}
\end{equation}

Let $u\in Z_k$ with $\norm{u} = r_k > 0$. We can choose uniformly an $\epsilon>0$ small enough such that
$\epsilon\abs{u}_{p^+}^{p^+} \leq \frac{1}{2p^+}\norm{u}^{p^-}$. Hence
$$\Phi(u) \geq \frac{1}{2p^+} \norm{u}^{p^-} - C\max\{ \abs{u}_{q(x)}^{q^-}, \abs{u}_{q(x)}^{q^+}\}.$$

If $\max\{\abs{u}_{q(x)}^{q^-}, \abs{u}_{q(x)}^{q^+}\} = \abs{u}_{q(x)}^{q^-}$, we choose $r_k=(2q^-C\alpha_k^{q^-}) ^ {\frac{1}{p^--q^-}}$, we get that
\begin{equation}
\begin{array}{rl}
\Phi(u) &\geq \frac{1}{2p^+} \norm{u}^{p^-} - C\abs{u}_{q(x)}^{q^-} \geq \frac{1}{2p^+} - C\alpha_k^{p^-}\norm{u}^{q^-}
    \\&\geq (\frac{1}{2p^+} - \frac{1}{2q^-})r_k^{p^-}.
\end{array}
\end{equation}
Since $q^- > p^+$ and $\alpha_k \to 0$, we obtain $b_k\to \infty$.

If $\max\{\abs{u}_{q(x)}^{q^-}, \abs{u}_{q(x)}^{q^+}\} = \abs{u}_{q(x)}^{q^+}$, we can similarly derive that $b_k\to \infty$. Hence we have shown $(B)$ holds. \qed
\vspace{12pt}

By far, we have shown the geometry conditions of the Fountain Theorem hold. In fact, in order to use
the Fountain Theorem to get our main result, we needn't verify the functional $\Phi$ satisfies the $(PS)_c$
condition for every $c>0$. It's enough if we could find a special $(PS)$ sequence for each $c_k$ and verify
the sequence we find has a convergence subsequence. Of course, the first step is to show the $(PS)_{c_k}$
sequence is bounded. Because there is no Ambrosetti-Rabinowits type condition, we couldn't give a direct proof.
Following the ideas in L. Jeanjean \cite{24} and W. Zou \cite{38}, we consider $\Phi$ as a member in a family
of functional. We will show almost all the functional in the family have bounded $(PS)$ sequences.
The following result due to W. Zou and M. Schechter \cite{38} is crucial for this purpose.

Let the notions be the same as in Theorem $\ref{theo2}$. Consider a
family of real $C^1$ functional $\Phi_\lambda$ of the form:
$\Phi_\lambda(u) := I(u)- \lambda J(u)$, where $\lambda \in \Lambda$
and $\Lambda$ is a compact interval in $[0, \infty)$. We make the
following assumptions:\\
 $(A_1) \Phi_\lambda$ maps bounded sets into
bounded sets uniformly for $\lambda \in \Lambda$. Moreover,
$\Phi_\lambda(-u) = \Phi_\lambda(u)$ for all $(\lambda,
u)\in\Lambda\times X$.\\
 $(A_2) J(u)\geq 0$  for all $u\in E;
I(u)\to\infty$ or $J(u)\to\infty$ as $\norm{u}\to\infty$.\\
 Let
\begin{equation}
a_k(\lambda) :=\max\{\Phi_\lambda(u) : u\in W_k, \norm{u}=\rho_k\},
\end{equation}

\begin{equation}
b_k(\lambda) :=\inf\{\Phi_\lambda(u): u\in Z_k, \norm{u}=r_k\}.
\end{equation}
Define
$$c_k(\lambda) = \inf_{\gamma\in \Gamma_k} \max_{u\in B_k} \Phi_\lambda (\gamma(u)),$$
$$\Gamma_k : =\{ \gamma\in C(B_k, X) : \gamma \text{ is odd and }\gamma|_{\partial B_k} = id \}.$$

\begin{theorem}\label{theo3}
Assume that $(A_1)$ and $(A_2)$ hold. If $b_k(\lambda) > a_k(\lambda)$ for all $\lambda \in \Lambda$,
then $c_k(\lambda) \geq b_k(\lambda)$ for all $\lambda\in\Lambda$. Moreover, for almost every $\lambda\in\Lambda$,
there exists a sequence of $\{u_n^k(\lambda)\}_{n=1}^\infty$ such that
$\sup_n\norm{u_n^k(\lambda)} < \infty, \Phi'_\lambda( u_n^k(\lambda))\to 0$ and $\Phi_\lambda(u_n^k(\lambda))\to c_k(\lambda)$ as $n\to\infty$.
\end{theorem}

Next, we let $I(u)=\int_{R^N}\frac{1}{\px} (\abs{Du}^{\px} + V(x)\abs{u}^{\px}) dx, J(u)=\int_{R^N} F(x, u)dx$
for $u\in E$ and $\Lambda=[1, 2]$. Under these terminologies, $\Phi(u)=\Phi_1(u)$. Under the assumptions of
Theorem $\ref{theo2}$. It's easy to see $(A_1)$ and $(A_2)$ hold.

\begin{lemma}\label{lemma38}
Under the assumptions of Theorem $\ref{theo2}$, $b_k(\lambda) > a_k(\lambda)$ for all $\lambda\in [1, 2]$ when $k$ is large enough.
\end{lemma}
{\bf Sketch of Proof.} Let $\rho_k > r_k >0$ large enough. Using
same reasoning, we can show that $a_k(\lambda) \leq 0$ and
$b_k(\lambda) \to \infty$ uniformly for $\lambda\in [1, 2]$ as
$k\to\infty$. Hence, we have shown the Lemma. Moreover,
$c_k(\lambda)\leq \sup_{u\in B_k} \Phi_\lambda(u)\leq \sup_{u\in
B_k}\Phi(u) = \max_{u\in B_k} \Phi_1(u) = \max_{u\in B_k}\Phi(u)
:=\ol{c_k} < \infty$. \qed\\

\begin{note}\label{note1}
Since $\Phi'_\lambda(u)$ is of type $(S_+)$ (Remark $\ref{remark1}$), we know that any bounded
$(PS)_{c(\lambda)}$ sequence of $\Phi_\lambda$ has a convergent subsequence which converges to a
critical point of $\Phi_\lambda$ with critical level $c(\lambda)$.
\end{note}

Now, using Theorem $\ref{theo3}$, we obtain that for almost every
$\lambda\in [1, 2]$, there exists a sequence of
$\{u_n^k(\lambda)\}_{n=1}^\infty$ such that
$\sup_n\norm{u_n^k(\lambda)} < \infty, \Phi'_\lambda(u_n^k(\lambda))
\to 0$ and $\Phi_\lambda(u_n^k(\lambda))\to c_k(\lambda)$ as
$n\to\infty$. Denote the set of these $\lambda$ by $\Lambda_0$. If
$1\in \Lambda_0$, we have found bounded $(PS)_{c_k}$ sequence for
the functional $\Phi$.

If $1\notin \Lambda_0$, we can choose a sequence $\{\lambda_n\} \subset \Lambda_0$
such that $\lambda_n\to 1$ decreasing. In view of Note $\ref{note1}$, for
each $\lambda\in\Lambda_0$, the bounded $(PS)_{c_k(\lambda)}$ sequence has a
convergent subsequence. We denote the limit by $u^k(\lambda)$. Accordingly,
$u^k(\lambda)$ is the critical point of the functional $\Phi_\lambda$ with
critical level $c_k(\lambda)$. Next, we are going to show the sequence
$\{ u^k(\lambda_n)\}_{n=1}^\infty$ is a bounded $(PS)_{c_k}$ sequence of
$\Phi$. For simplicity, we write $\{u^k(\lambda_n)\}$ as $\{u(\lambda_n)\}$.

In fact, we only need to show $\{ u(\lambda_n)\}$ is bounded. Indeed, if $\{u(\lambda)\}$ is bounded, we have
$$\Phi( u(\lambda_n) ) = \Phi_{\lambda_n}(u(\lambda_n)) + (1-\lambda_n) J(u(\lambda_n)) \to c_k,$$
$$\Phi'( u(\lambda_n) ) = \Phi'_{\lambda_n}(u(\lambda_n)) + (1-\lambda_n) J'(u(\lambda_n)) \to 0.$$
We have used the fact that $\Phi_\lambda, J$ map bounded sets into
bounded sets under the assumptions of Theorem 2.2.

\begin{lemma}
Under the assumption of Theorem 2.2, the aforementioned
$\{u(\lambda_n)\}$ is bounded.
\end{lemma}
{\bf Proof.} By contradiction. We assume $\norm{ u(\lambda_n) } \to \infty$ and consider $w_n = \frac{u(\lambda_n)}{ \norm{u(\lambda_n)} }$. Then up to a subsequence, we get that $w_n\rightharpoonup w$ in $E, w_n\to w$ in $L^{q(x)}(R^N)$ for $\px \leq q(x) \ll p^*(x), w_n\to w$ a.e. in $R^N$.

We first consider the case $w\neq 0$ in $E$. Since $\Phi'_{\lambda_n}( u(\lambda_n) ) = 0$, we have
$$\int_{R^N} \abs{Du(\lambda_n)}^{\px} + V(x)\abs{u(\lambda_n)}^{\px} dx = \lambda_n\int_{R^N} f(x, u(\lambda_n)) u(\lambda_n)dx.$$
Assume $\norm{u(\lambda_n)} > 1$. Dividing both sides by $\norm{u(\lambda_n)}^{p^+}$, we get
$$\int_{R^N} \frac{f(x, u(\lambda_n)) u(\lambda_n) }{ \norm{u(\lambda_n)}^{p^+} } dx \leq \frac{1}{\lambda_n} \leq 1.$$
Further, by Fatou's Lemma and $(C2)$, we have
$$\int_{R^N} \frac{f(x, u(\lambda_n)) u(\lambda_n) }{ \norm{u(\lambda_n)}^{p^+} } dx = \int_{R^N} \frac{f(x, u(\lambda_n)) u(\lambda_n)\abs{w_n(x)}^{p^+} }{ \norm{u_n(x)}^{p^+} } dx \to \infty,$$
a contradiction.

For the case $w=0$ in $E$, we define $\Phi_{\lambda_n}(t_nu(\lambda_n)) = \max_{t\in [0, 1]} \Phi_{\lambda_n} (tu(\lambda_n))$. Then for any $C>1, \ol{w_n} :=\frac{Cu(\lambda_n)}{ \norm{u(\lambda_n)} }$ and $n$ large enough, we have
$$
\begin{array}{rl}
\Phi_{\lambda_n}(t_nu(\lambda_n)) &\geq \Phi_{\lambda_n}( \ol{w_n})
\\& = \int_{R^N} \frac{1}{\px} (\abs{CDw_n}^{\px} + V(x)\abs{Cw_n}^{\px})dx-\lambda_n\int_{R^N}F(x, Cw_n)dx
\\& \geq\frac{1}{p^+}C^{p^-} - \lambda_n\int_{R^N} F(x, Cw_n)dx.
\end{array}
$$
Since $w_n\to 0$ a.e. in $R^N$ and $\lambda_n\in [1, 2]$, we have
$\lambda_n\int_{R^N} F(x, Cw_n)dx \to 0$ as $n\to\infty$. Since $C$
is arbitrary, we have $\Phi_{\lambda_n}(t_nu(\lambda_n)) \to\infty$
as $n\to\infty$. Consequently, we know $t_n\in (0, 1)$ when $n$ is
large enough, which implies
$\Phi'_{\lambda_n}(t_nu(\lambda_n))t_nu(\lambda_n) = 0$. Thus,
$$\Phi_{\lambda_n}(t_nu(\lambda_n)) - \frac{1}{p^-} \Phi'_{\lambda_n}(t_nu(\lambda_n))t_nu(\lambda_n) \to \infty,$$
which implies
\begin{align*}
\int_{R^N}(\frac{1}{\px} - \frac{1}{p^-}) ( \abs{t_nDu(\lambda_n)}^{\px} + V(x)\abs{t_nu(\lambda_n)}^{\px}) dx+
\\
 \lambda_n\int_{R^N} \frac{1}{p^-}f(x, t_nu(\lambda_n))t_nu(\lambda_n) - F(x, t_nu(\lambda_n))dx \to\infty.
 \end{align*}
So
$$\int_{R^N} \frac{1}{p^-} f(x, t_nu(\lambda_n))t_nu(\lambda_n) - F(x, t_nu(\lambda_n)) dx \to\infty.$$
However,
$$
\begin{array}{rl}
\Phi_{\lambda_n}(u(\lambda_n)) &=\Phi_{\lambda_n}(u(\lambda_n)) - \frac{1}{p^+}\Phi'_{\lambda_n}(u(\lambda_n))u(\lambda_n)
    \\&=\int_{R^N}(\frac{1}{\px} - \frac{1}{p^+})( \abs{Du(\lambda_n)}^{\px} + V(x)\abs{u(\lambda_n)}^{\px})dx
    \\&\ \ \  +\lambda_n\int_{R^N}\frac{1}{p^+}f(x, u(\lambda_n))u(\lambda_n) - F(x, u(\lambda_n))dx
    \\&\geq \lambda_n\int_{R^N}\frac{1}{p^+}f(x, u(\lambda_n))u(\lambda_n) - F(x, u(\lambda_n))dx.
\end{array}
$$
In view of $(C4)$, there exist two positive constants $C_1$ and
$C_2$ such that
$$
\begin{array}{rl}
\Phi_{\lambda_n}(u(\lambda_n)) &\geq \lambda_n\int_{R^N}\frac{1}{p^+}f(x, u(\lambda_n))u(\lambda_n) - F(x, u(\lambda_n))dx
        \\&\geq \lambda_n C_1\int_{R^N}\frac{1}{p^-}f(x, u(\lambda_n))u(\lambda_n) - F(x, u(\lambda_n))dx
        \\&\geq \lambda_n C_1C_2\int_{R^N} \frac{1}{p^-}f(x, t_nu(\lambda_n))t_nu(\lambda_n) - F(x, t_nu(\lambda_n))dx
        \\&\geq C\int_{R^N} \frac{1}{p^-}f(x, t_nu(\lambda_n))t_nu(\lambda_n) - F(x, t_nu(\lambda_n))dx\to\infty.
\end{array}
$$
However, for each $k$ large enough, $\Phi_{\lambda_n}(u(\lambda_n)) = c_k(\lambda_n) \leq \ol{c_k} < \infty$ (See Lemman $\ref{lemma38}$), a contradiction. \qed
\vspace{12pt}

{\bf Proof of Theorem $\ref{theo1}$.} By now, whether $1\in
\Lambda_0$ or not, we have found a special bounded $(PS)_{c_k}$
sequence $\{u^k(\lambda_n)\}_{n=1}^\infty$ for each $c_k$ in the
Fountain Theorem when $k$ is large enough. In view of Note
$\ref{note1}$, we know $\{u^k(\lambda_n)\}_{n=1}^\infty$ has a
convergent subsequence and $c_k$ is indeed an critical level of
$\Phi$ and Theorem 2.2 follows. \qed

\begin{remark}
We prove Theorem $\ref{theo1}$ in such a way because we want to
emphasize the procedure of finding critical points. First, we
consider the original functional and verify the functional satisfies
some geometry properties (e.g. Mountain Pass Geometry in \cite{24},
Fountain geometry in this paper, general linking geometry, etc) to
ensure prospective critical levels. Then, we consider our functional
as a member in a family of functionals. Some given structure
conditions on the family yield bounded $(PS)$ sequences for almost
all the functionals. Using the information supplied by these
functionals, we find special bounded $(PS)$ sequences for those
prospective critical levels. At last, we prove the special $(PS)$
sequences we find converge to critical points respectively up to
subsequences.
\end{remark}

\end{document}